\documentclass[10pt]{article}

\usepackage[english]{babel}
\usepackage[utf8x]{inputenc}
\usepackage[T1]{fontenc}
\usepackage[letterpaper,top=.8in,bottom=.8in,left=1in,right=1in,marginparwidth=1.5cm]{geometry}

\usepackage{amsmath,amsfonts,amsthm,amssymb,mathrsfs,dsfont} 
\usepackage{mathtools}
\usepackage{graphicx}
\usepackage[colorinlistoftodos]{todonotes}
\usepackage[colorlinks=true, allcolors=blue]{hyperref}
\usepackage[capitalize,nameinlink]{cleveref}
\usepackage{verbatim}
\usepackage{enumerate}

\global\long\def\N{\mathbb{N}}

\global\long\def\E{\mathbb{E}}

\global\long\def\Pr{\text{Pr}}

\newtheorem{theorem}{Theorem}[section]
\newtheorem*{namedtheorem}{\theoremname}
\newcommand{\theoremname}{testing}

\newtheorem{proposition}[theorem]{Proposition}

\newtheorem{corollary}[theorem]{Corollary}

\newtheorem{question}[theorem]{Question}

\newtheorem{conjecture}[theorem]{Conjecture}
\newtheorem*{question*}{Question}

\theoremstyle{definition}

\newtheorem{remark}[theorem]{Remark}

\theoremstyle{plain}

\usepackage{algpseudocode,algorithm,algorithmicx}

\allowdisplaybreaks

\title{Uniformity-independent minimum degree conditions for perfect matchings in hypergraphs}
\author{Asaf Ferber \thanks{Massachusetts Institute of Technology. Department of Mathematics. Email: {\tt ferbera@mit.edu}. Research is partially supported by an NSF grant 6935855.} \and Vishesh Jain\thanks{Massachusetts Institute of Technology. Department of Mathematics. Email: {\tt visheshj@mit.edu}}}
\date{}

\begin{document}
\maketitle
\begin{abstract} 
In this note, we prove that there exists a universal constant $c=\frac{43}{50}$ such that for every $k\in \mathbb{N}$ and every $d<k/2$, every $k$-uniform hypergraph on $n$ vertices and with minimum $d$-degree at least $(c+o_n(1))\binom{n-d}{k-d}$ contains a perfect matching. This is the first such bound which is independent of $k$, and therefore, improves all previously known bounds when $k$ is large. Our approach is based on combining the seminal work of Alon et al. with known bounds on a conjectured probabilistic inequality due to Feige. 
\end{abstract}
\section{Introduction}
In this note, we observe, following \cite{alon2012large}, the equivalence (\cref{prop:equivalence}) between the fractional, asymptotic version of the Erd\H{o}s matching conjecture in extremal hypergraph theory, and a probabilistic inequality that controls deviations from expectation of sums of non-negative i.i.d. random variables under only a first moment constraint. As a consequence, we are able to use existing work around this probabilistic inequality to significantly improve on the bounds of the extremal hypergraph theoretic problem (\cref{thm:main:matchings}) for certain choices of parameters, and to use existing work on the extremal hypergraph theoretic problem in other parameter ranges to strengthen the probabilistic inequality (\cref{corollary:improve-deviation-bound}). We will discuss these problems in more detail below.

\subsection{Perfect matchings in hypergraphs}
For $k\in\N$ with $k\geq2$, a \emph{$k$-uniform hypergraph} (or
\emph{$k$-graph} for short) $H$ is a pair $H=(V,E)$, where $V:=V(H)$
is a finite set of vertices, and $E:=E(H)\subseteq{V(H) \choose k}$
is a family of $k$-element subsets of $V$, referred to as the edges of
$H$. Given a $k$-graph $H$ and a subset $S\subseteq V$ with $0\leq|S|\leq k-1$,
we let $\deg_{H}(S)$ be the number of edges of $H$ containing
$S$. That is, 
$$\deg_H(S):=\left| \left\{e\in E(H) \mid S\subseteq e\right\}\right|.$$
The \emph{minimum $d$-degree} of $H$, $\delta_d(H)$, is defined as 
\[
\delta_{d}(H):=\min\left\{ \deg_{H}(S)\mid S\in{V(H) \choose d}\right\} .
\]
A collection of vertex disjoint edges of $H$ is called a \emph{matching}, and the number of edges in a matching is called the \emph{size}
of the matching. We say that a matching $M\subseteq E(H)$ is a \emph{perfect matching} if $|M|=|V|/k$ (in particular, this requires $|V|$ to be divisible
by $k$). 

For integers $n,k,d,s$ satisfying $0\leq d\leq k-1$ and $0\leq s\leq n/k$,
let $m_{d}^{s}(k,n)$ be the smallest integer $m$ such that every
$n$-vertex $k$-graph $H$ with $\delta_{d}(H)\geq m$ has a matching
of size $s$. Of particular interest is the case $s=n/k$ (that is, a minimum $d$-degree condition to enforce a perfect matching); for convenience, we will denote $m_{d}^{n/k}(k,n)$ by $m_{d}(k,n)$. 
The problem of determining $m_{d}(k,n)$ is a central theme in extremal graph theory and has attracted a lot of attention in the last few decades (see, e.g., the surveys \cite{rodl2010dirac, zhao2016recent} and the references therein). The case $m_1(2,n)$ goes back to a classical work of Dirac \cite{dirac1952some} from the 1950s, where he proved that every graph on $n$ vertices with minimum degree at least $n/2$ contains a Hamiltonian cycle (that is, a cycle of length $n$), and therefore, by keeping every other edge (assuming that $n$ is even), also contains a perfect matching. Ever since, problems of this type (that is, minimum degree conditions which guarantee existence of a perfect matching/Hamiltonian cycle in hypergraphs) are referred to as \emph{Dirac-type problems for hypergraphs}. One specific conjecture that has attracted a considerable amount of attention is the following: 
\begin{conjecture}
\label{conj:dirac-matching}
For all integers $1\leq d\leq k-1$ and $n$ which is divisible by $k$, 
\[
m_{d}(k,n)=\left(\max\left\{ \frac{1}{2},1-\left(1-\frac{1}{k}\right)^{k-d}\right\} +o_{n}(1)\right){n-d \choose k-d},
\]
where $o_{n}(1)$ stands for some function that tends to $0$ as $n$
tends to infinity. 
\end{conjecture}
\begin{remark}
Simple explicit constructions (see, e.g., \cite{alon2012large}) show that the right hand side is a lower bound on $m_{d}(k,n)$, so the content of the conjecture is that it is also an upper bound. 
\end{remark}
It is readily checked that for $d\geq k/2$, the maximum of the two
terms appearing in \cref{conj:dirac-matching} is $1/2$. In this case, \cref{conj:dirac-matching} is even known to hold in a stronger \emph{exact} form due to Treglown and Zhao \cite{treglown2012exact} (generalizing and
improving on previous work of R\"odl, Rucki\'nski, and Szemer\'edi \cite{rodl2006dirac} and
Pikhurko \cite{pikhurko2008perfect}); see also the work of Frankl and Kupavskii \cite{frankl1806erdHos}. On the other hand, for $d<k/2$, \cref{conj:dirac-matching} has been verified for
only a few cases (see the discussion in \cite{han2016perfect}) and, for example, even the case $d=1$ and $k=6$ is open. For $1\leq d<k/2$, the best known general
upper-bounds on $m_{d}(k,n)$ are the following (improving on earlier
results by H\`an, Person, and Schacht \cite{han2009perfect} and Markstr\"om and Rucki\'nski \cite{markstrom2011perfect}):
\begin{theorem}[K\"uhn, Osthus, and Townsend \cite{kuhn2014fractional}]  For integers $k\geq 3$, $1\leq d \leq k/2$ and $n\in k\N$,
\[
m_{d}(k,n)\leq\left(\frac{k-d}{k}-\frac{k-d-1}{k^{k-d}}+o_{n}(1)\right){n-d \choose k-d}.
\]
\end{theorem}
\begin{theorem}[Han, \cite{han2016perfect}]
\label{thm:han}
For integers $k\geq 3$, $1\leq d < k/2$ and $n\in k\N$,
\[
m_{d}(k,n)\leq\max\left\{ \delta(n,k,d)+1,(g(k,d)+o_{n}(1)){n-d \choose k-d}\right\} ,
\]
where $\delta(n,k,d) = (1/2 + o_n(1))\binom{n-d}{k-d}$ is some explicitly described function, and
\begin{align}
\label{eqn:gkd}
g(k,d):=1-\left(1-\frac{(k-d)(k-2d-1)}{(k-1)^{2}}\right)\left(1-\frac{1}{k}\right)^{k-d}.
\end{align}
\end{theorem}
For a detailed comparison of these bounds, we refer the reader to \cite{han2016perfect}. Here, we emphasize that if we think of $k$ as large and $d = o(k)$, then both the above bounds only show that 
$$m_{d}(k,n) \leq (1-o_k(1))\binom{n-d}{k-d},$$
whereas \cref{conj:dirac-matching} predicts that
$$m_{d}(k,n) \leq\left(1-e^{-1} + o_k(1)\right) \binom{n-d}{k-d}.$$

In this note, we take a first step towards \cref{conj:dirac-matching} by showing that the above statement holds with a non-zero absolute constant (unfortunately, smaller than $e^{-1}$), thereby substantially improving the above bounds in the particularly interesting case when $d$ is small compared to $k$.
\begin{theorem}
\label{thm:main:matchings}
 For integers $k\geq 3$, $1\leq d \leq k/2$ and $n$ which is divisible by $k$,
$$m_{d}(k,n) \leq \left(1 - \frac{7}{50} + o_n(1)\right)\binom{n-d}{k-d}.$$
\end{theorem}
\begin{remark}
Although our proof of this theorem does not involve any significant new ideas -- indeed, after this note first appeared on the arXiv, we were informed by Andrey Kupavskii that such an improvement using the same methods is also implicit in his work \cite{frankl1806erdHos} with Peter Frankl on the Erd\H{o}s Matching Conjecture -- we believe that the statement itself is important enough to be worth isolating and highlighting. It was also pointed out to us by Andrey Kupavskii that the equivalence between the problems established in \cref{prop:equivalence} already appears (with the same proof) in the work of {\L}uczak, Mieczkowska, and {\v{S}}ileikis \cite{luczak2017maximal}, although in this work, no mention is made either of Feige's conjecture, or the implication for Dirac-type thresholds for perfect matchings in hypergraphs.   
\end{remark}
For the remainder of this note, we define
$$m_d(k):= \limsup_{n\in k\N, n\to \infty}\frac{m_{d}(k,n)}{\binom{n-d}{k-d}}.$$

\subsection{Fractional perfect matchings in hypergraphs}

A \emph{fractional matching} in a $k$-graph $H=(V,E)$ is a function
$w:E\to[0,1]$ such that for every $v\in V$, $\sum_{e\ni v}w(e)\leq1$ (observe that if $w:E\to \{0,1\}$, then the same condition gives a matching).
The \emph{size} of a fractional matching is defined to be $\sum_{e\in E}w(e)$.
We say that $w$ is a \emph{perfect fractional matching} if its size
is $|V|/k$ (or equivalently, if $\sum_{e\ni v}w(e)=1$ for all $v\in V$). Note that as every matching is also a fractional matching (but not vice versa), it follows that the size of the largest fractional matching in $H$ is atEven though our proof of this theorem does not involve any significant new ideas, we believe that the least as large as the size of the largest matching in $H$. 

Analogously to the (integer) matching case, for an integer $0\leq d\leq k-1$ and a real number $0\leq s\leq n/k$, we let $f_{d}^{s}(k,n)$ denote the
smallest integer $m$ such that every $n$-vertex $k$-graph $H$
with $\delta_{d}(H)\geq m$ has a fractional matching of size
$s$. As before, we will denote $f_{d}^{n/k}(k,n)$ simply by $f_{d}(k,n)$.
In their seminal work, Alon, Frankl, Huang, R\"odl, Ruci\'nski, and Sudakov \cite{alon2012large} showed that in order to prove \cref{conj:dirac-matching}, it suffices to prove the corresponding variant for fractional perfect matchings. More precisely, setting
$$f_{d}(k):= \limsup_{n\to \infty}\frac{f_{d}(k,n)}{\binom{n-d}{k-d}},$$
they proved the following: 
\begin{theorem}[Theorem 1.1 in \cite{alon2012large}]
Fix integers $k,d$ (with $k\geq3$ and $1\leq d < k/2$) and a real $\alpha>0$.
Then, there exists $n_{0}:=n_{0}(k,d,\alpha)\in\N$ such that for
all $n\geq n_{0}$, all $k$-graphs $H$ on $n$ vertices with $\delta_{d}(H)\geq (f_{d}(k) + \alpha)\binom{n-d}{k-d}$ contain a perfect matching. 
\end{theorem}
In particular, it follows that for the above range of parameters
\begin{align}
\label{eqn:1}
m_{d}(k) = f_{d}(k)+o(1).
\end{align}
Hence, it suffices to study the following (presumably easier) analogous conjecture for fractional perfect matchings. 
\begin{conjecture}
\label{conjecture:fpm}
For integers $k\geq 3$ and $1\leq d \leq k-1$,
$$f_d(k) = 1-\left(1-\frac{1}{k}\right)^{k-d}.$$
\end{conjecture}
Once again, the lower bound on $f_d(k)$ follows by a simple, explicit construction. As a means to prove the above conjecture, the authors in \cite{alon2012large} provided a further reduction to the `$d=0$' case.
\begin{proposition}[Proposition 1.1 in \cite{alon2012large}] 
\label{prop:reduce-to-em}
For integers $k\geq 3$, $1\leq d \leq k-1$, and $n\geq k$,
$$f_d(k,n) \leq f_0^{n/k}(k-d,n-d).$$
\end{proposition}

In particular, this shows that 
it suffices to prove the following, which is a fractional version of the classical matching conjecture of Erd\H{o}s \cite{erdHos1965problem} (see \cite{alon2012large} for a more general statement and further discussion). 
\begin{conjecture}
\label{conjecture:fractional-erdos-matching}
For integers $\ell, d\geq 1$ and $s := (m+d)/(\ell + d)$,
$$\limsup_{m\to \infty}\frac{f_{0}^{s}(\ell,m)}{\binom{m}{\ell}} \leq 1-\left(1-\frac{1}{\ell+d}\right)^{\ell}.$$
\end{conjecture}
Henceforth, we will denote the quantity appearing on the left hand side of the above inequality by $f^{d}(\ell)$; this should \emph{not} be confused with $f_{d}(\ell)$, which was defined earlier. With this notation, \cref{prop:reduce-to-em} shows that for any $k\geq 3$ and $1\leq d \leq k-1$, 
\begin{align}
\label{eqn:2}
f_{d}(k) \leq f^{d}(k-d).
\end{align}

\subsection{Small deviation inequalities for sums of independent random variables}
Using an elegant LP duality argument (which we will essentially replicate in the next section), the authors in \cite{alon2012large} showed that \cref{conjecture:fractional-erdos-matching} follows from a conjectured probabilistic inequality due to Samuels \cite{samuels1966chebyshev}. Since Samuels proved his conjecture in some special cases, a corresponding resolution of \cref{conjecture:fractional-erdos-matching} in a few cases was obtained in \cite{alon2012large}. Here, based on the approach from \cite{alon2012large}, we will show that \cref{conjecture:fractional-erdos-matching} is, in fact, equivalent to special cases of a conjectured probabilistic inequality (of a similar nature to Samuels' conjecture) due to Feige  \cite{feige2006sums}. Given this, our main result will follow by using known bounds on Feige's conjecture, as we now discuss.\\

In \cite{feige2006sums}, Feige asked the following question: 
\begin{question} Suppose that $X_1,\dots,X_n$ are  independent, non-negative random variables, each of which has mean $1$. How large can 
$\Pr[X_1+\dots+X_n \geq n+1] $ 
be?\end{question} 

Note that Markov's inequality gives the upper bound $n/(n+1) = 1-1/(n+1)$ which is essentially useless. Indeed, Feige conjectured that a much better bound should hold. 
\begin{conjecture}[Feige, \cite{feige2006sums}]
\label{conjecture:feige}
Let $X_1,\dots,X_n$ be independent, non-negative random variables, each of which has mean $1$. Then, for all $d\geq 1$,
$$\Pr[X_1+\dots+X_n \geq n+d] \leq 1- \left(1-\frac{1}{n+d}\right)^{n}.$$
\end{conjecture}
Observe that the conjectured bound is attained when $X_1,\dots,X_n$ are i.i.d. random variables 
which take on the value $n+d$ with probability $1/(n+d)$, and the value $0$ otherwise. 
In the same paper \cite{feige2006sums}, Feige also proved that \cref{conjecture:feige} holds with $1- (1-1/(n+d))^{n}$ replaced by $12/13$, and this was later improved to $7/8$ by He, Zhang, and Zhang \cite{he2010bounding}. To the best of our knowledge, the current best bound in this direction is $43/50$ due to Garnett \cite{garnett2018small}. We remark that \cite{feige2006sums, he2010bounding} appeared much before \cite{alon2012large}.

For finding perfect matchings in hypergraphs, the only case in Feige's inequality which is of interest to us is when $X_1,\dots,X_n$ are i.i.d. random variables. Accordingly, for $\ell \in \N$ and $d > 0$, we define the following quantity:
$$\Theta^{d}(\ell):= \sup\Pr[X_1 + \dots + X_\ell \geq \ell + d],$$
where the supremum is taken over all collections of non-negative i.i.d. random variables $X_1,\dots,X_\ell$ with mean $1$. In particular, it follows from Garnett's bound that for all $\ell \in \N$ and $d\geq 1$, 
\begin{align}
\label{eqn:3}
\Theta^{d}(\ell) \leq \frac{43}{50}.
\end{align}
\begin{remark}
\label{rmk:limit}
For later use, we note that $\liminf_{\epsilon \downarrow 0}\Theta^{d-\epsilon}(\ell) = \Theta^{d}(\ell)$. Indeed, the direction $\geq $ is immediate, whereas for the direction $\leq $, we use the following observation: for any $\delta > 0$, and any non-negative mean $1$ random variable $X$, let $Y_{\delta}$ denote the random variable which takes on the value $0$ with probability $\delta$, and is distributed as $(1-\delta)^{-1}\cdot X$ otherwise. Then, $Y_{\delta}$ is a non-negative mean $1$ random variable, and letting $Y_1,\dots,Y_\ell$ (resp. $X_1,\dots,X_\ell$) denote i.i.d. copies of $Y_\delta$ (resp. $X$),
$$\Pr[Y_1 + \dots + Y_\ell \geq \ell+d] \geq (1-\delta)^{\ell}\Pr[X_1+\dots + X_\ell \geq  (\ell + d)(1-\delta)]\geq (1-\delta)^{\ell} \liminf_{\epsilon \downarrow 0}\Theta^{d-\epsilon}(\ell),$$
so that the desired inequality follows by taking $\delta \downarrow 0$.
\end{remark}
In the next section, we will prove the following.
\begin{theorem}
\label{prop:equivalence}
For any integers $\ell,d \geq 1$ we have
$$f^{d}(\ell) = \Theta^{d}(\ell).$$
\end{theorem}
It is clear that this proposition implies \cref{thm:main:matchings}. Indeed, for integers $k\geq 3$ and $1\leq d \leq k/2$, we have 
$$m_{d}(k) = f_{d}(k) \leq f^{d}(k-d) = \Theta^{d}(k-d) \leq \frac{43}{50},$$
where the first equality is \cref{eqn:1}, the second inequality is \cref{eqn:2}, the third equality is \cref{prop:equivalence}, and the last inequality is \cref{eqn:3}.\\

On the other hand, we can use \cref{prop:equivalence} along with (the proof of) \cref{thm:han} to obtain improved bounds for the i.i.d. version of Feige's inequality for sufficiently large deviations. It is shown in \cite{han2016perfect} (see the last computation in the proof of Theorem 1.5 there) that for integers $n,k$ with $n > 3k \geq 6$,
$$f_{0}^{n/k}(k-d,n-d) \leq \left(g(k,d)+o_{n}(1)\right)\binom{n-d}{k-d}$$
(where $g(k,d)$ is the function defined in \cref{eqn:gkd}), from which it immediately follows that
\begin{align}
    f^{d}(\ell) \leq g(\ell + d,d).
\end{align}
Hence, we obtain:
\begin{corollary}
\label{corollary:improve-deviation-bound}
Let $X_1,\dots,X_\ell$ be i.i.d. non-negative random variables each of which has mean $1$. Then, for all integers $d\geq 1$,
$$\Pr[X_1+\dots+X_\ell \geq \ell+d]\leq g(\ell+d,d).$$
\end{corollary}
Note that this bound does indeed improve on Markov's inequality/Feige's inequality/Hoeffding's inequality for $d$ which is sufficiently large compared to $\ell$. For instance, one may directly check (as is done in Corollary 1.4 in \cite{han2016perfect}) that for $0.73\ell \leq d < \ell$, $g(\ell + d,d)< 1/2$, whereas the above mentioned inequalities give a bound which is worse than $1/2$. 
\section{Proof of \cref{prop:equivalence}}
We begin by showing that $f^{d}(\ell) \leq \Theta^{d}(\ell)$. The proof is exactly the same as the proof of Theorem 2.1 in \cite{alon2012large}, with the only change being replacing the application of Samuels' conjecture in the last step with Feige's conjecture. Before proceeding with the details we need to introduce some notation. For an $\ell$-graph $H$, we will denote the size of the largest fractional matching by $\nu^{*}(H)$. Then, by LP-duality, we have
$$\nu^*(H) = \tau^*(H),$$
where $\tau^*(H)$ denotes the size of the smallest \emph{fractional vertex cover} in $H$. Recall that a fractional vertex cover of $H = (V,E)$ is a function $t\colon V\to [0,1]$ such that for every $e\in E$, we have $\sum_{v\in e}t(v)\geq 1$. The \emph{weight} of a fractional vertex cover is defined to be $\sum_{v\in V}t(v)$. 

Let $H$ be an $\ell$-graph on a vertex set $V$ of size $m$, and
suppose that $\nu^{*}(H)= xm$, where 
$$0 < x \leq \frac{1+ \frac{d}{m}}{\ell + d}.$$ 
Then, since $\tau^{*}(H)= xm$
there exists a weight function $t:V\to[0,1]$ such that $\sum_{v\in V}t(v)\leq xm$
and for every edge $e\in E(H)$ we have $\sum_{v\in e}t(v)\geq1$. Let $v_{1},\dots,v_{\ell}$ be a sequence of vertices of $H$, chosen
independently and uniformly from $V$. For each $i\in[\ell]$, we
define a random variable $X_{i}:=t(v_{i})$. Hence, $X_{1},\dots,X_{\ell}$
are i.i.d. random variables with mean 
\[
\mu:=\frac{1}{m}\sum_{v\in V}t(v)=x.
\]
We will now use bounds on the deviation of the sum $X_1 +\dots + X_\ell$ to bound the number of edges in $H$. 

Since $\sum_{e\ni v}t(v)\geq1$ for all $e\in E$, the number $N$
of $\ell$-element subsets $S$ of $V$ with $\sum_{v\in S}t(v)\geq 1$
is an upper bound on the number of edges of $H$. Let $N_{1}$
denote the number of all $\ell$-element sequences of vertices of
$V$ whose sum of weights is at least $1$, and let $N_{2}$
denote the number of all $\ell$-element sequences of distinct vertices
of $V$ whose sum of weights is at least $1$. Note that
$N_{1}-N_{2}\leq{\ell \choose 2}m^{\ell-1}=O_{\ell}(m^{\ell-1})$
(since $N_{1}-N_{2}$ is at most the number of $\ell$-element sequences
in which at least one vertex appears twice) and $N_{2}=\ell!N$. Therefore,
\[
\Pr\left[\sum_{i=1}^{\ell}t(v_{i})\geq 1\right]=\frac{N_{1}}{m^{\ell}}=\frac{N_{2}+O_{\ell}(m^{\ell-1})}{{m \choose \ell}\ell!}=(1+o_{m}(1))\frac{N}{{m \choose \ell}}.
\]
Moreover, we also have
\begin{align*}
\Pr\left[\sum_{i=1}^{\ell}t(v_{i}) \geq 1\right]
&=\Pr\left[\sum_{i=1}^{\ell}X_{i} \geq 1\right]\\
&= \Pr\left[ \sum_{i=1}^{\ell} x^{-1}X_i \geq x^{-1}\right]\\
&\leq \Pr\left[\sum_{i=1}^{\ell} x^{-1}X_i \geq \ell + d-o_{m}(1)\right]\\
&\leq \Theta^{d - o_m(1)}(\ell).
\end{align*}
The desired conclusion now follows from \cref{rmk:limit} by taking the liminf of the right hand side as $m \to \infty$.\\ 

We will now prove the reverse inequality $\Theta^{d}(\ell) \leq f^{d}(\ell)$, which follows by `reversing' the above proof. Let $X_1,\dots,X_\ell$ be non-negative i.i.d. random variables with mean $1$. For the purpose of bounding $\Theta^{d}(\ell)$ from above, it clearly suffices to assume that each $X_i$ is supported in $[0,\ell +d]$. By a standard approximation argument, we may further restrict our attention only to those distributions which are supported on finitely many points, and which take on each value with probability equal to some rational number. Hence, $X_1,\dots, X_\ell \sim X$, where $X$ is supported on $x_1,\dots,x_a \in [0,\ell +d]$, and $\Pr[X = x_i] = b_i/m'$, for some $m',b_i\in \N$ with $\sum_{i=1}^{a}b_i = m'$.

We now construct an $\ell$-graph $H=(V,E)$ as follows. Let $V=[m]$, where $m=rm'$ for $r\in \N$. Let $t\colon V\to [0,1]$ be defined by sending the first $rb_1$ elements of $[m]$ to $x_1/(\ell + d)$, the next $rb_2$ elements of $[m]$ to $x_2/(\ell + d)$ and so on, and define 
$$E:=\left\{S \in \binom{|V|}{\ell} : \sum_{v\in S}t(v) \geq 1\right\}.$$
Then, by construction, $t$ is a fractional vertex cover of $H$, so that 
$$\nu^*(H) = \tau^*(H) \leq \sum_{v\in V}t(v) = \frac{\E[X]}{\ell+d}\cdot m = \frac{1}{\ell+d} m.$$

As before, let $v_1,\dots,v_\ell$ be a sequence of vertices of $H$, chosen independently and uniformly from $V$, and for each $i\in[\ell]$, let $X_i := t(v_i)$. Then, $X_1,\dots,X_\ell$ are i.i.d. random variables, and observe from the definition of $t$ that $(\ell+d)\cdot X_i \sim X$. Moreover, by the same argument as above, it follows that 
\begin{align*}
\Pr\left[\sum_{i=1}^{\ell} (\ell+d)\cdot X_i \geq \ell + d\right]
&= \Pr\left[\sum_{i=1}^{\ell}t(v_i) \geq 1 \right]\\
&= (1+o_m(1))\frac{|E|}{\binom{m}{\ell}}\\
&\leq (1+o_m(1))\frac{f_{0}^{(m+d)/(\ell+d)}(\ell,m)}{\binom{m}{\ell}}.
\end{align*}
Finally, the desired conclusion follows by noting that the liminf of the right hand side (as $r\to \infty$) is at most $f^{d}(\ell)$.
\\

\noindent {\bf Acknowledgement.} We would like to thank Mathias Schacht for communicating this problem to us, and Andrey Kupavskii for bringing references \cite{frankl1806erdHos,luczak2017maximal} to our attention.  
\bibliographystyle{abbrv}
\bibliography{main}
\end{document}